\documentclass[letterpaper, 10 pt, conference]{ieeeconf}  

\IEEEoverridecommandlockouts                              

\usepackage{graphicx,algorithm,algorithmic,amsthm,amsmath,amsfonts,verbatim,mathtools,enumerate,bm,hyperref,stmaryrd,cite}
\usepackage{tikz, pgfplots, subcaption}
\usepackage{siunitx}
\usepackage{booktabs}
\usepackage{multirow}
\pgfplotsset{compat=1.18}
\usepgfplotslibrary{fillbetween}

\allowdisplaybreaks

\newcommand{\norm}[1]{\left\lVert#1\right\rVert}
\newcommand{\mnorm}[1]{{\left\vert\kern-0.25ex\left\vert\kern-0.25ex\left\vert #1 
    \right\vert\kern-0.25ex\right\vert\kern-0.25ex\right\vert}}

\title{\LARGE \bf Nonlinear Trajectory Optimization Models for Energy-Sharing UAV-UGV Systems with Multiple Task Locations
}

\author{Minsen Yuan, Amanuel Adane, James Humann, and Yue Yu
\thanks{M. Yuan and Y. Yu are with the Department of Aerospace Engineering and Mechanics, University of Minnesota, Minneapolis, MN 55455, USA ({\tt\small \{yuan0450, yuey\}@umn.edu}). A. Adane is with the Duffield College of Engineering, Cornell University, Ithaca, NY 14853, USA ({\tt\small aa2888@cornell.edu}). J. Humann is with DEVCOM Army Research Laboratory ({\tt\small james.d.humann.civ@army.mil}). Y. Yu would like to thank Samet Uzun for helpful early discussions.}
}

\begin{document}

\maketitle
\thispagestyle{empty}
\pagestyle{empty}

\begin{abstract}

Energy-sharing UAV-UGV systems extend the endurance of Uncrewed Aerial Vehicles (UAVs) by leveraging Uncrewed Ground Vehicles (UGVs) as mobile charging stations, enabling persistent autonomy in infrastructure-sparse environments. Trajectory optimization for these systems is often challenging due to UGVs' terrain access constraints and the discrete nature of task scheduling. We propose a smooth nonlinear program model for the joint trajectory optimization for these systems. Unlike existing models, the proposed model allows smooth parameterization of UGVs' terrain access constraints and supports partial UAV recharging. Further, it introduces a smooth approximation of disjunctive constraints that eliminates the need for computationally expensive integer programming and enables efficient solutions via nonlinear programming algorithms. We demonstrate the proposed model on a one-UAV-one-UGV system with multiple task locations. Compared with mixed-integer nonlinear programs, this model reduces the computation time by orders of magnitude.


\end{abstract}

\definecolor{NLP_color}{rgb}{1,0,0}
\definecolor{MINLP_color}{rgb}{0.301, 0.745, 0.933}

\definecolor{AL_Lp_color}{rgb}{1,0,0}
\definecolor{AL_Logsumexp_color}{rgb}{0.301, 0.745, 0.933}
\definecolor{IPM_Lp_color}{rgb}{0, 1, 0}
\definecolor{IPM_Logsumexp_color}{rgb}{0.9290,0.6940,0.1250}

\section{Introduction}
Collaborative UAV-UGV systems provide a powerful platform to integrate the heterogeneous strengths of aerial and ground robots via \emph{energy-sharing}. In these systems, Uncrewed Ground Vehicles (UGVs) serve as mobile docking and charging stations that effectively extended the endurance of Uncrewed Aerial Vehicles (UAVs) via \emph{sharing} energy stored on UGVs. By combining UGVs' superior payload capacity together with UAVs' agility and elevated perspective, these systems enable a wide range of robotic applications in remote and challenging environments where access to traditional energy infrastructure is sparse. Examples of these applications include mapping, disaster response, surveillance, and precision agriculture \cite{ding2021review,chai2024cooperative,munasinghe2024comprehensive}. 

Coordinating the trajectories for energy-sharing UAV-UGV systems presents a unique challenge, especially in missions with multiple tasks locations. This difficulty is twofold. First, The UGV's terrain access constraints restrict the UGVs to a known road network. Second, scheduling the the order in which the UAVs and UGVs completes their individual tasks often require optimizing discrete decisions. When these discrete decisions are coupled with the terrain access constraints, the resulting optimization is often nonconvex and computationally expensive to solve. 

One approach to optimize UAV-UGV trajectories is via the \emph{two-echelon routing problem}, where the trajectories of the UAV and UGV are computed in two distinct phases. One approach is to first determining the UAV trajectory, then optimize the UGV trajectory to support the UAV  \cite{mathew2015}. Another approach is to first determine the UGV trajectory that provide a set of feasible rendezvous locations to support UAV recharging, then plan the UAV trajectory based these locations \cite{maini2015, mondal2025cooperative}. Early results uses a greedy algorithm to generate the rendezvous locations, then compute the UAV trajectory accordingly \cite{maini2015,maini2019}. More recent results augment this approach with an asynchronous team framework \cite{ramasamy2024} and task allocation heuristics based on minimum set covering \cite{mondal2025cooperative}. Another approach is to decouple the planning of UAV and UGV trajectory using reachable sets constraints \cite{kim2024decoupled}. In all such cases, the resulting trajectory quality depends on how to partition the the planning between the UAV and UGV.

To directly address the coupling between UAVs and UGVs, an alternative approach is to jointly optimize their trajectories. On approach for such optimization is mixed integer linear programs (MILP), which computes an optimal UAV-UGV trajectory via optimization over integer variables \cite{wang2019,manyam2019}. Since MILP is often computationally expensive, there have been several strategies to simplify the problem under stronger assumptions. For example, if the potential rendezvous locations are known, genetic algorithm provides an efficient alternative to MILP \cite{eker2025}. If the order of UAV and UGV actions are known, then optimizing UAV-UGV trajectories reduces to a convex second-order cone program \cite{diller2026}. A more recent direction is deep reinforcement learning methods to train encoder-decoder based transformer networks that generates the UAV-UGV trajectories \cite{wu2023,mondal2024,li2025}. These trajectories can not only minimize mission time but also account for stochastic battery usage by constraining on the risk of mid-mission battery depletion \cite{mondal2025}.   

There are several limitations in the existing models for UAV–UGV trajectory optimization. First, the UGV road network is often modeled as a graph with discretized nodes, causing the planning complexity to grow rapidly with the number of nodes \cite{mathew2015,maini2015,maini2019,wang2019,manyam2019,mondal2025cooperative,ramasamy2024,eker2025}. Second, most existing models assume that the UAV  battery is fully recharged before takeoff \cite{maini2019,diller2026,mondal2024,mondal2025cooperative,ramasamy2024,eker2025}. While this assumption is compatible with battery replacement, it is not suitable for wireless charging that allow partial recharge. Finally, MILP methods rely on integer variables to model discrete decisions, leading to solution times that grow exponentially with the number of integer variables \cite{wang2019,manyam2019,maini2019,mondal2025cooperative,mathew2015,maini2015,ramasamy2024}.

To address these limitations, we propose a novel \emph{nonlinear} trajectory optimization model for energy-sharing UAV–UGV systems with multiple task locations. We base this model on three key features. First, we model the UGV road network as the union of multiple continuous paths joined at one junction point. This approach avoids the need for node discretization along each path and allows smooth parameterization of UGV's terrain access constraints. Second, we model the dynamics of UAV battery via distinct charging and discharging curves. Compared with traditional models that always only allow full recharges, this approach allows  partial recharging and increase the flexibility of UAV trajectories. Third, we develop a smoothing approach to model disjunctive constraints that are necessary for scheduling optimal task completion and battery recharging. We first transform the disjunctive constraints into nonsmooth constraints defined by pointwise minimum function, then propose a smooth approximation of these constraints based on \(\ell_p\)-norms. This approximation enables efficient solution with smooth nonlinear optimization methods, bypassing the need for computationally expensive integer programming. We demonstrate this model on a one-UAV–one-UGV system with multiple task locations. Compared with mixed-integer nonlinear programming, this model reduces solution time from hours to minutes while maintaining robust success rates in numerical experiments.

\section{Trajectory Optimization Models for Energy-Sharing UAV-UGV Systems}
\label{sec: opt model}
We model the trajectory optimization for energy-sharing UAV-UGV system---where the UGV serve as mobile charging stations for the UAV---as a constrained optimization problem. This model includes the optimization variables used to parameterize the UAV-UGV trajectory, the objective function, and the physical and operational constraints of the UAV-UGV system. 
\subsection{Trajectory Variables}
We consider three classes of variables when optimizing the UAV-UGV trajectory: variables shared by both trajectories, variables specific to the UAV trajectory, and variables specific to the UGV trajectory. To simplify the notation, we let \(\llbracket k \rrbracket \coloneqq \{1, 2, \ldots, k\}\) for any \(k \in \mathbb{N}\).
\paragraph{Shared Variables} Let \(N\in\mathbb{N}\) denote the total number of time stamps along the UAV trajectory, which is the same as the total number of time steps along the UGV trajectory. We let \(s_k\in\mathbb{R}_{\geq 0}\) denote the time duration between the \(k\)-th and the \((k+1)\)-th time stamp for all \(k\in \llbracket N\rrbracket\).  
\paragraph{UAV Variables} We parameterize the UAV trajectory using its planar position and its remaining time of flight supported by the current battery level. We let \(r_k^{\texttt{A}}\in\mathbb{R}^2\) denote the projection of UAV's 3D position onto the xy-plane at the \(k\)-th time stamp for all \(k\in\llbracket N\rrbracket\), and \(e_k\in\mathbb{R}_{\geq 0}\) denote the remaining time of flight of the UAV at time step \(k\).
\paragraph{UGV Variables} We assume that the UGV only moves on a road network in the shape of a star graph, consisting of \( m^{\texttt{G}} \) arms that intersect only at a single junction point. We parameterize the UGV position trajectory using two variables. 
We let \( r_k^{\texttt{G}} \in \mathbb{R}^2 \) denote the planar position of the UGV at the 
\(k\)-th time stamp. In addition, we let 
\( p_k \in \mathbb{R}_{\ge 0}^{m^{\texttt{G}}} \) 
denote the UGV’s position on the star graph with \( m^{\texttt{G}} \) arms. 
Each entry of \( p_k \) corresponds to one arm of the star graph, and the value of a positive entry is the distance the UGV traveled from the graph center along that arm to its current position. At any time stamp, \( p_k \) has at most one positive entry, indicating that the UGV occupies at most one arm of the graph at a time.
\subsection{Objective Function}
We choose the total mission time as the objective function. Minimizing total time directly captures the operational efficiency of the UAV–UGV system, encouraging timely completion of all required tasks while implicitly balancing travel, waiting, and charging time. Since the time durations between consecutive time stamps are optimization variables, the total mission time is given by
\begin{equation}
   \textstyle  \sum_{k=1}^{N-1} s_k.
\end{equation}

\subsection{Trajectory Constraints}
We consider three classes of constraints for the joint UAV–UGV trajectory: constraints that depend only on UAV variables, constraints that depend only on UGV variables, and constraints that couple both UAV and UGV variables.
\paragraph{UAV Constraints}
We let \(\overline{r}_0\in\mathbb{R}^2\) and \(\overline{r}_f\in\mathbb{R}^2\) denote the initial and final position of the UAV trajectory, respectively. We consider the following initial and final constraints
\begin{equation}
    r_1^{\texttt{A}}=\overline{r}_0, \enskip r_N^{\texttt{A}}=\overline{r}_f.
\end{equation}
In addition, we consider the following constraint on the speed of the UAV
\begin{equation}
    \norm{r_{k+1}^{\texttt{A}} - r_k^{\texttt{A}}}_2 \leq v_{\texttt{max}}^{\texttt{A}} s_k,
\end{equation}
for all \( k \in\llbracket N-1\rrbracket \), where \(v_{\texttt{max}}^{\texttt{A}}\in\mathbb{R}_{>0}\) is the maximum UAV speed. Note that the UAV position is three-dimensional, whereas the above constraint accounts only for planar motion. This is because we assume that the UAV operates at a constant altitude, and that the time required for vertical motion (\emph{e.g.}, during takeoff and landing) is negligible compared with the time spent in horizontal motion.

We consider the case where the UAV must visit a set of task locations along its trajectory (e.g., to monitor areas of interest). Let \( m^{\texttt{A}} \in \mathbb{N} \) denote the total number of task locations, and let 
\( a_1, \ldots, a_{m^{\texttt{A}}} \in \mathbb{R}^2 \) denote the locations of these task locations on the \(xy\)-plane. We consider the following logical constraints:
\begin{equation}
    \exists\, k \in\llbracket N\rrbracket  \text{ s.t. } r_k^{\texttt{A}} = a_i,
\end{equation}
for all \( i \in\llbracket m^{\texttt{A}}\rrbracket  \). 
These constraints ensure that the UAV visits each task point at least once along its trajectory. 

\paragraph{UGV Constraints}
We consider the case where the initial and final position of the UGV coincide with those of the UAV, namely,
\begin{equation}
    r_1^{\texttt{G}}=\overline{r}_0, \enskip r_N^{\texttt{G}}=\overline{r}_f.
\end{equation}
The UGV can only move along a star-shaped road network. At the \(k\)-th time stamp, \( r_k^{\texttt{G}} \in \mathbb{R}^2 \) denotes the planar position 
of the UGV, and \( p_k \in \mathbb{R}_{\ge 0}^{m^{\texttt{G}}} \) denotes its position 
along the arms of the star graph. 
In particular, \( p_k = \mathbf{0}_{m^{\texttt{G}}} \) indicates that the UGV is located 
at the junction point (i.e., the center of the star graph). 
Moreover, the UGV lies on the \(j\)-th arm if and only if \( [p_k]_j > 0 \).
We introduce a nonlinear mapping \( g : \mathbb{R}_{\ge 0}^{m^{\texttt{G}}} \to \mathbb{R}^2 \) such that
\begin{equation}
    r_k^{\texttt{G}} = g(p_k).
\end{equation}
Let \( p_{\texttt{max}} \in \mathbb{R}_{\ge 0}^{m^{\texttt{G}}} \) denote a parameter vector whose 
\(j\)-th entry represents the maximum allowable distance of the UGV from the junction point along 
the \(j\)-th arm, corresponding to the furthest task point on that arm. We consider the following constraints on \( p_k \):
\begin{equation}\label{eqn: UGV terrain}
\begin{aligned}
    & \mathbf{0}_{m^{\texttt{G}}} \leq p_k \leq p_{\texttt{max}},  \enskip r_k^{\texttt{G}} = g(p_k),\\
    & p_k^\top 
    \bigl( \mathbf{1}_{m^{\texttt{G}}}\mathbf{1}_{m^{\texttt{G}}}^\top - I_{m^{\texttt{G}}} \bigr)
    p_k = 0,
\end{aligned}
\end{equation}
for all \(k \in \llbracket N \rrbracket \). These constraints ensure that the UGV occupies at most one arm of the star graph at any time stamp, and that its planar position \(r_k^{\texttt{G}}\) is coupled with the vector \( p_k \) through the nonlinear 
mapping \(g\).

Furthermore, we consider the following constraints on the speed of the UGV
\begin{equation}\label{eqn: UGV speed}
    \norm{p_{k+1}-p_k}_1\leq v_{\texttt{max}}^{\texttt{G}} s_k
\end{equation}
for all \(k\in\llbracket N-1\rrbracket\), where \(v_{\texttt{max}}^{\texttt{G}}\in\mathbb{R}_{>0}\) is the maximum UGV speed. These constraint in \eqref{eqn: UGV speed}, when combined with the constraint in \eqref{eqn: UGV terrain}, ensure that the distance traveled by the UGV is upper bound by the product of its maximum speed and the time used. 

Similar to the UAV, the UGV also must visit a set of task locations. Since the UGV can only move on a road network modeled as a star graph, it is sufficient to specify the furthest point along each arm that the UGV must reach. Without loss of generality, we take the endpoint of each arm to be the corresponding furthest point. To this end, we consider the following constraints:
\begin{equation}
    \exists\, k \in \llbracket N \rrbracket  \text{ s.t. } 
    b_j^\top ( p_k - p_{\texttt{max}} ) = 0,
\end{equation}
for all \( j \in \llbracket m^{\texttt{G}} \rrbracket \), where 
\( b_j \in \mathbb{R}^{m^{\texttt{G}}} \) denotes the \(j\)-th column of the identity matrix 
\( I_{m^{\texttt{G}}} \). These constraints ensure that the UGV trajectory reaches the end of each arm of the road network at least once along its trajectory.

\begin{figure}[!t]
\centering
\begin{tikzpicture}

\node[anchor=south west, inner sep=0] (img) 
    {\includegraphics[width=0.65\linewidth]{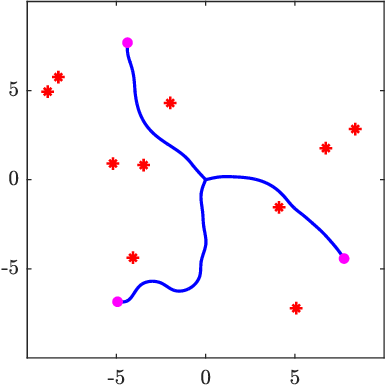}};

\begin{scope}[x={(img.south east)}, y={(img.north west)}]

\node at (0.5,-0.05) {\footnotesize $ X\;(\mathrm{km})$};
\node[rotate=90] at (-0.05,0.5) {\footnotesize $ Y\;(\mathrm{km})$};

\node at (0.28,0.90) {\footnotesize $\overline{r}_0$};
\node at (0.25,0.20) {\footnotesize $\overline{r}_f$};

\end{scope}

\end{tikzpicture}
\caption{Example problem with $m^{\texttt{A}}=10$ UAV task locations (red stars) and $m^{\texttt{G}}=3$ arms, whose associated UGV task locations are shown as magenta circles.}
\label{fig: an_example}
\end{figure}
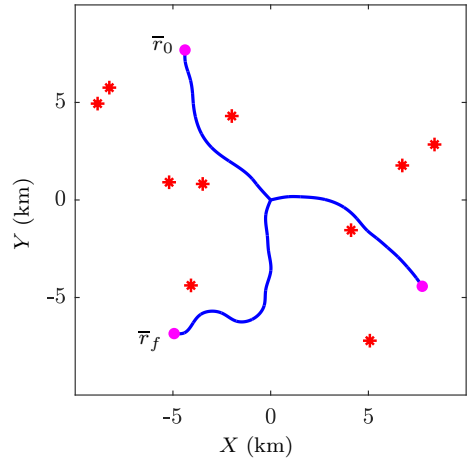

\paragraph{Coupling Constraints}
We consider the following bounds on the time duration between two consecutive time stamps, which constrain both UAV and UGV trajectories:
\begin{equation}
    s_{\texttt{min}} \le s_k \le s_{\texttt{max}},
\end{equation}
for all \( k \in \llbracket N-1 \rrbracket \), where 
\( s_{\texttt{min}} \in \mathbb{R}_{>0} \) and 
\( s_{\texttt{max}} \in \mathbb{R}_{>0} \) denote the lower and upper bounds on the time duration between consecutive time stamps.

Furthermore, the UAV trajectory and UGV trajectory are coupled through the UAV battery dynamics. We first consider the following bounds on the UAV battery level (measured in remaining flight time)
\begin{equation}\label{eqn: energy bounds}
    e_{\texttt{min}}\leq e_k\leq e_{\texttt{max}},
\end{equation}
for all \( k \in \llbracket N\rrbracket \), where \(e_{\texttt{min}}\in\mathbb{R}_{\geq 0}\) and \(e_{\texttt{max}}\in\mathbb{R}_{>0}\) denote the minimum and maximum allowed values for the battery level. In addition, the UAV battery dynamics include two modes: charging and discharging. The UAV can discharge its battery at any location, but it can recharge only when it lands on the UGV, in which case their positions coincide. We model these requirements using the following constraints:
\begin{equation}
    \label{eqn: battery_dynamics}
    \left\{
    \begin{aligned}
         e_{k+1} & = \min( e_k + \kappa s_k, e_{\texttt{max}} ), \\
         r_k^{\texttt{A}} & = r_k^{\texttt{G}}, \\
         r_{k+1}^{\texttt{A}} & = r_{k+1}^{\texttt{G}}
    \end{aligned}
    \right\}
    \lor
    \left\{
        e_{k+1} = e_k - s_k
    \right\},
\end{equation}
for all \( k \in \llbracket N-1 \rrbracket \), where 
\( \kappa \in \mathbb{R}_{>0} \) denotes the UAV battery charging rate when the UAV is landed on the UGV. These constraints ensure that, between two consecutive points along the UAV trajectory, the UAV battery level either decreases, increases, or reaches capacity. In the latter case, the UAV and UGV positions coincide at both trajectory points. Moreover, matching positions at two consecutive trajectory points imply that the UAV and UGV positions coincide for all times between the corresponding time stamps. 

In practice, we find that the constraint 
\( e_{k+1} = \min\!\left( e_k + \kappa s_k,\; e_{\texttt{max}} \right) \) 
is often overly restrictive for numerical solvers, particularly when combined with the bounding constraints in~\eqref{eqn: energy bounds}. As a remedy, we relax this constraint as
\begin{equation}\label{eqn: relaxed UAV charging}
    e_{k+1} \leq \min\left( e_k + \kappa s_k, e_{\texttt{max}} \right).
\end{equation}
Intuitively, this relaxation allows the battery level to be lower than the ideal charging outcome, hence providing additional flexibility for the solver without compromising the practical feasibility of the resulting solution.

\section{Nonlinear Smoothing for Disjunctive Constraints}
The challenge in solving the trajectory optimization problem proposed in Section~\ref{sec: opt model} comes from the constraints. We can divided the constraints discussed in Section~\ref{sec: opt model} into two groups. The first group consists of constraints defined by smooth functions:
\begin{equation}\label{eqn: smooth constr}
    \begin{aligned}
&r_1^{\texttt{A}}=r_1^{\texttt{G}}=\overline{r}_0, \enskip r_N^{\texttt{A}}=r_N^{\texttt{G}}=\overline{r}_f,\enskip e_1=e_{\texttt{max}},\\
&\norm{r_{k+1}^{\texttt{A}} - r_k^{\texttt{A}}}_2 \leq v_{\texttt{max}}^{\texttt{A}} s_k,\enskip  k \in\llbracket N-1\rrbracket,\\
&\norm{p_{k+1}-p_k}_1\leq v_{\texttt{max}}^{\texttt{G}} s_k,\enskip  k \in\llbracket N-1\rrbracket,\\
&  e_{\texttt{min}}\leq e_{k+1}\leq e_{\texttt{max}},\enskip s_{\texttt{min}} \leq s_k \le s_{\texttt{max}},\enskip k \in\llbracket N-1\rrbracket,\\
& r_k^{\texttt{G}} = g(p_k), \enskip \mathbf{0}_{m^{\texttt{G}}} \leq p_k \leq p_{\texttt{max}},\enskip k \in\llbracket N\rrbracket,\\
& p_k^\top 
    \bigl( \mathbf{1}_{m^{\texttt{G}}}\mathbf{1}_{m^{\texttt{G}}}^\top - I_{m^{\texttt{G}}} \bigr)
    p_k = 0,\enskip k \in\llbracket N\rrbracket.
    \end{aligned}
\end{equation}
The constraints in \eqref{eqn: smooth constr} are compatible with many algorithms for smooth nonlinear programs, such as interior point methods and augmented Lagrangian methods. 

The second group consists of disjunctive constraints, given as follows:
\begin{equation}\label{eqn: disjunctive constr}
    \begin{aligned}
        & \exists\, k \in\llbracket N\rrbracket  \text{ s.t. } r_k^{\texttt{A}} = a_i,\enskip i \in\llbracket m^{\texttt{A}} \rrbracket,\\
        & \exists\, k \in \llbracket N \rrbracket  \text{ s.t. } 
    b_j^\top ( p_k - p_{\texttt{max}} ) = 0, \enskip j\in \llbracket m^{\texttt{G}} \rrbracket,\\
    &\left\{
    \begin{aligned}
         e_{k+1} & \leq e_k + \kappa s_k \\
         r_k^{\texttt{A}} & = r_k^{\texttt{G}} \\
         r_{k+1}^{\texttt{A}} & = r_{k+1}^{\texttt{G}}
    \end{aligned}
    \right\}
    \lor
    \left\{
        e_{k+1} = e_k - s_k
    \right\},\, k\in \llbracket N-1 \rrbracket.
    \end{aligned}
\end{equation}
Notice that here we changed the constraint in \eqref{eqn: relaxed UAV charging} to \(e_{k+1}  \leq  e_k + \kappa s_k\). This change is lossless when we impose the constraints in \eqref{eqn: smooth constr} and \eqref{eqn: disjunctive constr} together since the constraints in \eqref{eqn: smooth constr} already ensures that \(e_{k+1}\leq e_{\texttt{max}}\) for all \(k\in\llbracket N\rrbracket\).

The disjunctive constraints in~\eqref{eqn: disjunctive constr} pose unique challenges for optimization, as they rely on logical OR operations and induce a disconnected feasible solution set. We first discuss how to model these constraints using discrete variables, leading to a mixed-integer nonlinear programming approach. We then propose an alternative approach that first reformulates the constraints in~\eqref{eqn: disjunctive constr} as nonsmooth constraints and subsequently approximates them with smooth nonlinear functions that are compatible with nonlinear programming.
  
\subsection{Disjunctive Constraints via Discrete Variables}
A classical approach to model the disjunctive constraints in \eqref{eqn: disjunctive constr} using discrete binary variables. In particular, we can reformulate these constraints as follows:
\begin{equation}\label{eqn: MINLP constr}
\begin{aligned}
    &  \norm{r_k^{\texttt{A}} - a_i}_\infty \leq \mu (1-U_{ki}),\enskip k\in \llbracket N \rrbracket, i\in \llbracket m^{\texttt{A}} \rrbracket,\\
    & |b_j^\top ( p_k - p_{\texttt{max}})| \leq \mu (1-V_{kj}),\enskip k\in \llbracket N \rrbracket, j\in \llbracket m^{\texttt{G}} \rrbracket,\\
    & |e_{k+1}^{\texttt{A}}-e_k^{\texttt{A}}+s_k|\leq \mu (1-W_k),\enskip k\in \llbracket N-1 \rrbracket,\\
    & e_{k+1}^{\texttt{A}}-e_k^{\texttt{A}}-\kappa s_k\leq \mu W_k,\enskip k\in \llbracket N-1 \rrbracket,\\
    & \norm{\begin{matrix}
              r_k^{\texttt{A}}-r_k^{\texttt{G}}\\
              r_{k+1}^{\texttt{A}}-r_{k+1}^{\texttt{G}}
          \end{matrix}}_\infty\leq \mu W_k,\enskip  k\in \llbracket N-1 \rrbracket,\\
     & \textstyle \sum_{k=1}^N U_{ki}=1, U_{ki}\in\{0, 1\}, \enskip  k\in \llbracket N \rrbracket, i\in \llbracket m^{\texttt{A}} \rrbracket,\\
     & \textstyle \sum_{k=1}^N V_{kj}=1, V_{kj}\in\{0, 1\}, \enskip  k\in \llbracket N \rrbracket, j\in \llbracket m^{\texttt{G}} \rrbracket,\\
     & W_k\in\{0, 1\}, \enskip  k\in \llbracket N-1 \rrbracket,\\
\end{aligned}
\end{equation}
where \(\mu\gg 1\) is a large positive scalar. The idea is to first express each disjunctive constraint as a finite set of candidate conditions, one of which must be satisfied. Next, we introduce binary variables to encode the selection of these candidate conditions, and ensure feasibility by forcing the sum of the associated binary variables to equal one.

We formulate the UAV-UGV trajectory optimization problem a \emph{mixed integer nonlinear program (MINLP)}. This program contains the following variables
\begin{equation}\label{eqn: MINLP var}
\left\{r_k^{\texttt{A}}, r_k^{\texttt{G}},
        e_k,  p_k, 
        \{U_{k i}\}_{i=1}^{m^{\texttt{A}}},
        \{V_{k j}\}_{j=1}^{m^{\texttt{G}}}
    \right\}_{k=1}^N \cup \left\{
        s_k, W_k
    \right\}_{k=1}^{N-1}.
\end{equation}
We formulate this MINLP as follows
\begin{equation}\label{opt: MINLP}
\begin{array}{ll}
    \underset{\text{ Variables in \eqref{eqn: MINLP var}}}{\text{minimize}}
    &  \sum_{k=1}^{N-1} s_k \\[4pt]
    \;\; \text{subject to}
    & \text{constraints in \eqref{eqn: smooth constr} and \eqref{eqn: MINLP constr}}.
\end{array}
\end{equation}
We can solve MINLP above using branch-and-bound–based methods combined with nonlinear programming algorithms. For details on models, algorithms, and practical solution methods for MINLP, we refer interested readers to~\cite{grossmann1997mixed}. 

\subsection{Disjunctive Constraints via Nonlinear Smoothing}
One limitation of the discrete-variable approach is that it leads to an exponential growth in the number of possible values for binary variables, which often makes scalable real-time solutions impractical. 
As an alternative, we introduce a continuous modeling approach for the disjunctive constraints in~\eqref{eqn: disjunctive constr} that avoids the use of discrete variables. To this end, we first reformulate the constraints in~\eqref{eqn: disjunctive constr} as follows: 
\begin{equation}\label{eqn: nonsmooth constr}
\begin{aligned}  
  &\underset{k \in\llbracket N\rrbracket}{\min} \norm{r_k^{\texttt{A}} - a_i}_2=0,\enskip i \in\llbracket m^{\texttt{A}} \rrbracket, \\
    & \underset{k \in\llbracket N\rrbracket}{\min}
    |b_j^\top ( p_k - p_{\texttt{max}})| = 0, \enskip j\in \llbracket m^{\texttt{G}} \rrbracket, \\
&\min\left(|e_{k+1}^{\texttt{A}}-e_k^{\texttt{A}}+s_k|,\norm{\begin{matrix}
              \sigma_\delta(e_{k+1}^{\texttt{A}}-e_k^{\texttt{A}}-\kappa s_k)\\
              r_k^{\texttt{A}}-r_k^{\texttt{G}}\\
              r_{k+1}^{\texttt{A}}-r_{k+1}^{\texttt{G}}
          \end{matrix}}_2
          \right)= 0,\\
         & k\in \llbracket N-1 \rrbracket,    
\end{aligned}
\end{equation}
where 
\begin{equation}
    \textstyle \sigma_\delta(\alpha) = \begin{cases}
        0, & \alpha\leq 0,\\
        \frac{1}{2}\alpha^2, & 0\leq \alpha\leq \delta,\\
        \delta \alpha -\frac{1}{2}\delta^2, & \alpha>\delta
    \end{cases}
\end{equation}
and \(\delta\in\mathbb{R}_{>0}\) is a parameter with small positive value. The key idea is to reformulate each disjunctive constraint using pointwise minimum of the violation of candidate conditions. Here, the function \(\psi_\delta\) provides a smooth measure of the violation of inequality constraint.

Note that the constraints in \eqref{eqn: nonsmooth constr} are not compatible with algorithms for smooth nonlinear programs, due to the nonsmooth pointwise minimum function. A common approach to approximate the pointwise minimum function is via the \emph{log-sum-exp} function. Given \(c_1, c_2, \ldots, c_n\in\mathbb{R}\), the approximation is as follows
\begin{equation}
    \underset{k\in\llbracket  n\rrbracket}{\min} \, c_k  \textstyle\approx -\frac{1}{\tau}\ln \left(\sum_{i=1}^n \exp(-\tau c_k) \right),\label{eqn: log-sum-exp}
\end{equation}
where \(\tau\in\mathbb{R}_{>0}\) is a positive scaling parameter. The idea of this approximation is that, as \(\tau\) increases, the exponential term corresponds to the smallest entry will dominate the sum and cancel with the logarithm function.

However, the log-sum-exp function often causes numerical instabilities due to the rapid growth of the exponential function. To over come this limitation, we propose to replace the log-sum-exp function a smooth function constructed based on the \(\ell_p\)-norm. In particular, let \(c\) denote the vector whose \(k\)-th entry is \(c_k\), we let \(c^{-1}\) denote the elementwise reciprocal of \(c\), assuming \(c_k>0\) for all \(k\). We can show that
\begin{equation}\label{eqn: lp limit}
\lim\limits_{p\to \infty, p\in\mathbb{N}} \sqrt[2p]{n}/\norm{c^{-1}}_{2p} =1/\norm{c^{-1}}_\infty=\underset{k\in\llbracket n\rrbracket}{\min}\, c_k. 
\end{equation}
Here the coefficient \(\sqrt[2p]{n}\) normalize the magnitude of the norm \(\norm{c^{-1}}_{2p}\). This limit implies that, for a sufficiently large integer \(p\), the \(\ell_p\)-norm provides a smooth formula to approximate the pointwise minimum function. To avoid the singularity case where \(c\) contains a zero entry, we propose the following approximation formula
\begin{equation}\label{eqn: lp smoothing}
\underset{k\in\llbracket  n\rrbracket}{\min} \, c_k  \approx \textstyle \left(\frac{1}{n}\sum_{i=1}^n (c_k^2+\epsilon^2)^{-p}\right)^{-\frac{1}{2p}}-\epsilon,  
\end{equation}
where \(\epsilon\in\mathbb{R}_{>0}\) is a small positive parameter to ensure the function is well defined, and \(p\in\mathbb{N}\) is a reasonably large integer (in practice, \(p=3\) provides a satisfying approximation in simulation). Notice that \eqref{eqn: lp smoothing} reduces to \eqref{eqn: lp limit} if \(\epsilon=0\). Similar approximation has been used for logical specifications in optimization \cite{uzun2024optimization,uzun2025motion}.


With either one of the formulas above, we can approximate the nonsmooth constraints in \eqref{eqn: nonsmooth constr} as smooth constraints and optimize the UAV-UGV trajectory optimization by solving a smooth nonlinear optimization problem. This \emph{nonlinear program (NLP)} contains the following variables: 
\begin{equation}\label{eqn: NLP var}
\left\{
        r_k^{\texttt{A}}, r_k^{\texttt{G}},
        e_k,  p_k
    \right\}_{k=1}^N \cup \left\{
        s_k
    \right\}_{k=1}^{N-1}.
\end{equation}
We formulate this NLP as follows
\begin{equation}\label{opt: NLP}
\begin{array}{ll}
\underset{\text{Variables in \eqref{eqn: NLP var}}}{\text{minimize}}
    &  \sum_{k=1}^{N-1} s_k \\
    \;\; \text{subject to}
    & \text{constraints in \eqref{eqn: smooth constr} and smooth approx. } \\
    & \text{of the constraints in \eqref{eqn: nonsmooth constr} via \eqref{eqn: log-sum-exp} or  \eqref{eqn: lp smoothing}.}
\end{array}
\end{equation}
Note that all of the variables in optimization~\eqref{opt: NLP} are continuous (\emph{i.e.}, no discrete-valued variables) and all of the functions that appear in optimization~\eqref{opt: NLP} are differentiable. As a result, in principle one can solve this optimization using standard algorithms for nonlinear programs (NLP). However, the smooth nonlinear functions in \eqref{eqn: lp smoothing} are often ill-conditioned. In practice, we observe that the augmented Lagrangian method provide consistent robust performance for solving ill-conditioned NLP \cite{chen2020convergence}.

\section{Numerical Simulation}


We demonstrate the proposed NLP model on a one-UAV-one-UGV system with multiple UAV and UGV task locations. Fig.~\ref{fig: an_example} shows the route map. Within the NLP framework, we compare the performance of different smoothing functions and nonlinear programming algorithms. We also compare the NLP model with the MINLP model to illustrate its scalability.

\subsection{Problem Setup}
We consider a problem with $m^{\texttt{G}}=3$, and the UAV and UGV share the same initial ($\overline{r}_0$) and final ($\overline{r}_f$) positions, as illustrated in Fig.~\ref{fig: an_example}. 
We consider $e_{\max}=0.4\,\mathrm{h}$, $v_{\max}^{\texttt{A}}=36\,\mathrm{km/h}$, $v_{\max}^{\texttt{G}}=16.2\,\mathrm{km/h}$, and $\kappa=1.5$. 
Additionally, we set $e_{\min}=0$, $s_{\min}=0$, and choose a sufficiently large upper bound $s_{\max}=10\,\mathrm{h}$ for the inequality constraints in~\eqref{eqn: smooth constr}.
We perform all simulations on the Minnesota Supercomputing Institute cluster (\url{https://www.msi.umn.edu/}), which uses AMD EPYC 7702 processors. Each simulation uses a single CPU core with 8~GB of allocated memory.

\subsection{Comparison of Different Algorithms for Nonlinear Programs}\label{subsec: NLP algorithm}
We compare the performance of the proposed NLP model in~\eqref{opt: NLP} using two commonly used algorithms for NLP: the Augmented Lagrangian Method (ALM) ~\cite{doi:10.1137/0312021,chen2020convergence} and the Interior-Point Method (IPM)~\cite{wright1997primal}. For the ALM, we use the \texttt{L-BFGS} method (implemented by \texttt{fminunc} in \texttt{MATLAB}) to minimize the augmented Lagrangian (see \cite{doi:10.1137/0312021}). For the IPM, we use the \texttt{fmincon} function in \texttt{MATLAB} with limited-memory Hessian approximation. Note that the convergence of NLP algorithms is sensitive to initialization~\cite{betts1998survey,yuan2025filtering}. This issue is particularly relevant in optimization~\eqref{opt: NLP}, since approximating nonsmooth functions with smooth ones leads to almost discontinuous gradients.
To mitigate this issue, we warm-start the NLP (ALM and IPM) by initializing the UAV's and UGV's position trajectories by setting their rendezvous locations as the projections of the UAV task locations onto the UGV's road network.

Fig.~\ref{fig: Convergence_Obj_Con_ten_aoi} shows the convergence of the ALM and IPM under the $\ell_p$-norm approximation in~\eqref{eqn: lp smoothing} and the \emph{log-sum-exp} approximation in~\eqref{eqn: log-sum-exp}. We compare these methods using 100 problem instances with randomly generated UAV task locations ($m^{\texttt{A}}=10$). The constraint violation refers to the sum of the violations of all constraints in~\eqref{eqn: smooth constr} and~\eqref{eqn: disjunctive constr}. We set \(\delta = 1\) in~\eqref{eqn: nonsmooth constr}, \(\epsilon = 10^{-3}\) and \(p = 3\) in~\eqref{eqn: lp smoothing}, and \(\tau = 10^{2}\) in~\eqref{eqn: log-sum-exp}. 
Among all the methods compared, the combination of ALM with the $\ell_p$-norm approximation achieves the best overall performance. With the same smoothing functions, ALM consistently converges to better local optima than IPM in terms of both objective value and constraint violation. When combined with the ALM, the $\ell_p$-norm approximation achieves lower constraint violation while maintaining a comparable objective value compared with the \emph{log-sum-exp} approximation.

\begin{figure}[!htp]
\centering

\begin{tikzpicture}[baseline]
\pgfplotslegendfromname{sharedlegend}
\end{tikzpicture}

\begin{subfigure}{\columnwidth}
\centering
\pgfplotsset{every tick label/.append style={font=\scriptsize}}
\begin{tikzpicture}
\begin{semilogyaxis}[
    xlabel = \footnotesize Time (s), ylabel=\footnotesize Objective Function Value,
    width= 0.9\linewidth,
    height = 0.45\linewidth,
    xmin = 0, xmax = 100,
    log ticks with fixed point,
    ytick = {3, 6, 9, 12},
    legend to name = sharedlegend, 
    legend style={
        legend columns=2,
        at={(0.5,-0.15)},
        anchor=north west,
        /tikz/column 2/.style={
            column sep=10pt},
        font=\footnotesize
    },
    legend image post style={line width=1.5pt} 
]

\definecolor{AL_Lp_color}{rgb}{1,0,0}
\addplot [mark=none, color=AL_Lp_color,line width=1pt] table[x=averageTimes,y=objMedianVals] {AL_Lp_ten_aoi.dat};
\addlegendentry{ALM w/ \(\ell_p\)-norm}

\definecolor{AL_Logsumexp_color}{rgb}{0.301, 0.745, 0.933}
\addplot [mark=none, color=AL_Logsumexp_color,line width=1pt] table[x=averageTimes,y=objMedianVals] {AL_Logsumexp_ten_aoi.dat};
\addlegendentry{ALM w/ \emph{log-sum-exp}}

\definecolor{IPM_Lp_color}{rgb}{0, 1, 0}
\addplot [mark=none, color=IPM_Lp_color,line width=1pt] table[x=averageTimes,y=objMedianVals] {IPM_Lp_ten_aoi.dat};
\addlegendentry{IPM w/ \(\ell_p\)-norm}

\definecolor{IPM_Logsumexp_color}{rgb}{0.9290,0.6940,0.1250}
\addplot [mark=none, color=IPM_Logsumexp_color,line width=1pt] table[x=averageTimes,y=objMedianVals] {IPM_Logsumexp_ten_aoi.dat};
\addlegendentry{IPM w/ \emph{log-sum-exp}}



\addplot[
    only marks,
    mark=*,
    mark size=0.8pt,
    color=IPM_Logsumexp_color,
    error bars/.cd,
        y dir=both,
        y explicit,
        error bar style={line width=0.1pt},
        error mark=-,  
        error mark options={
            line width=0.8pt,
            mark size=1.5pt,
            rotate=90
        },
]
table[
    x=averageTimes,
    y=objMedianVals,
    y error plus expr=\thisrow{objUpperMedianVals}-\thisrow{objMedianVals},
    y error minus expr=\thisrow{objMedianVals}-\thisrow{objLowerMedianVals},
]{IPM_Logsumexp_ten_aoi.dat};


\addplot[
    only marks,
    mark=*,
    mark size=0.8pt,
    color=IPM_Lp_color,
    error bars/.cd,
        y dir=both,
        y explicit,
        error bar style={line width=0.1pt},
        error mark=-,  
        error mark options={
            line width=0.8pt,
            mark size=1.5pt,
            rotate=90
        },
]
table[
    x=averageTimes,
    y=objMedianVals,
    y error plus expr=\thisrow{objUpperMedianVals}-\thisrow{objMedianVals},
    y error minus expr=\thisrow{objMedianVals}-\thisrow{objLowerMedianVals},
]{IPM_Lp_ten_aoi.dat};


\addplot[
    only marks,
    mark=*,
    mark size=0.8pt,
    color=AL_Logsumexp_color,
    error bars/.cd,
        y dir=both,
        y explicit,
        error bar style={line width=0.1pt},
        error mark=-,  
        error mark options={
            line width=0.8pt,
            mark size=1.5pt,
            rotate=90
        },
]
table[
    x=averageTimes,
    y=objMedianVals,
    y error plus expr=\thisrow{objUpperMedianVals}-\thisrow{objMedianVals},
    y error minus expr=\thisrow{objMedianVals}-\thisrow{objLowerMedianVals},
]{AL_Logsumexp_ten_aoi.dat};



\addplot[
    only marks,
    mark=*,
    mark size=0.8pt,
    color=AL_Lp_color,
    error bars/.cd,
        y dir=both,
        y explicit,
        error bar style={line width=0.1pt},
        error mark=-,  
        error mark options={
            line width=0.8pt,
            mark size=1.5pt,
            rotate=90
        },
]
table[
    x=averageTimes,
    y=objMedianVals,
    y error plus expr=\thisrow{objUpperMedianVals}-\thisrow{objMedianVals},
    y error minus expr=\thisrow{objMedianVals}-\thisrow{objLowerMedianVals},
]{AL_Lp_ten_aoi.dat};


\end{semilogyaxis}
\end{tikzpicture}
\caption{Convergence of objective function value.}
\end{subfigure}

\begin{subfigure}{\columnwidth}
\centering
\pgfplotsset{every tick label/.append style={font=\scriptsize}}
\begin{tikzpicture}
\begin{semilogyaxis}[
    xlabel = \footnotesize Time (s), ylabel=\footnotesize Constraint Violation,
    width= 0.9\linewidth,
    height = 0.45\linewidth,
    xmin = 1e-2, xmax = 100,
    ytick = {1e-6,  1e-4, 1e-2, 1e0,  1e2}
]

\definecolor{IPM_Logsumexp_color}{rgb}{0.9290,0.6940,0.1250}
\addplot [mark=none, color=IPM_Logsumexp_color,line width=1pt] table[x=averageTimes,y=conMedianVals] {IPM_Logsumexp_ten_aoi.dat};


\addplot[
    only marks,
    mark=*,
    mark size=0.8pt,
    color=IPM_Logsumexp_color,
    error bars/.cd,
        y dir=both,
        y explicit,
        error bar style={line width=0.1pt},
        error mark=-,  
        error mark options={
            line width=0.8pt,
            mark size=1.5pt,
            rotate=90
        },
]
table[
    x=averageTimes,
    y=conMedianVals,
    y error plus expr=\thisrow{conUpperMedianVals}-\thisrow{conMedianVals},
    y error minus expr=\thisrow{conMedianVals}-\thisrow{conLowerMedianVals},
]{IPM_Logsumexp_ten_aoi.dat};


\definecolor{IPM_Lp_color}{rgb}{0, 1, 0}
\addplot [mark=none, color=IPM_Lp_color,line width=1pt] table[x=averageTimes,y=conMedianVals] {IPM_Lp_ten_aoi.dat};

\addplot[
    only marks,
    mark=*,
    mark size=0.8pt,
    color=IPM_Lp_color,
    error bars/.cd,
        y dir=both,
        y explicit,
        error bar style={line width=0.1pt},
        error mark=-,  
        error mark options={
            line width=0.8pt,
            mark size=1.5pt,
            rotate=90
        },
]
table[
    x=averageTimes,
    y=conMedianVals,
    y error plus expr=\thisrow{conUpperMedianVals}-\thisrow{conMedianVals},
    y error minus expr=\thisrow{conMedianVals}-\thisrow{conLowerMedianVals},
]{IPM_Lp_ten_aoi.dat};


\definecolor{AL_Logsumexp_color}{rgb}{0.301, 0.745, 0.933}
\addplot [mark=none, color=AL_Logsumexp_color,line width=1pt] table[x=averageTimes,y=conMedianVals] {AL_Logsumexp_ten_aoi.dat};

\addplot[
    only marks,
    mark=*,
    mark size=0.8pt,
    color=AL_Logsumexp_color,
    error bars/.cd,
        y dir=both,
        y explicit,
        error bar style={line width=0.1pt},
        error mark=-,  
        error mark options={
            line width=0.8pt,
            mark size=1.5pt,
            rotate=90
        },
]
table[
    x=averageTimes,
    y=conMedianVals,
    y error plus expr=\thisrow{conUpperMedianVals}-\thisrow{conMedianVals},
    y error minus expr=\thisrow{conMedianVals}-\thisrow{conLowerMedianVals},
]{AL_Logsumexp_ten_aoi.dat};


\definecolor{AL_Lp_color}{rgb}{1,0,0}
\addplot [mark=none, color=AL_Lp_color,line width=1pt] table[x=averageTimes,y=conMedianVals] {AL_Lp_ten_aoi.dat};

\addplot[
    only marks,
    mark=*,
    mark size=0.8pt,
    color=AL_Lp_color,
    error bars/.cd,
        y dir=both,
        y explicit,
        error bar style={line width=0.1pt},
        error mark=-,  
        error mark options={
            line width=0.8pt,
            mark size=1.5pt,
            rotate=90
        },
]
table[
    x=averageTimes,
    y=conMedianVals,
    y error plus expr=\thisrow{conUpperMedianVals}-\thisrow{conMedianVals},
    y error minus expr=\thisrow{conMedianVals}-\thisrow{conLowerMedianVals},
]{AL_Lp_ten_aoi.dat};


\end{semilogyaxis}
\end{tikzpicture}
\caption{Convergence of constraint violation.}
\end{subfigure}

\caption{Convergence of different NLP algorithms and smoothing functions over 100 problem instances with randomly generated  UAV task locations with \(m^\texttt{A} = 10\). The solid lines represent the median values of the simulations, while the lower and upper bars indicate the interquartile range, spanning from the 0.25 quantile to the 0.75 quantile.}
\label{fig: Convergence_Obj_Con_ten_aoi}
\end{figure}
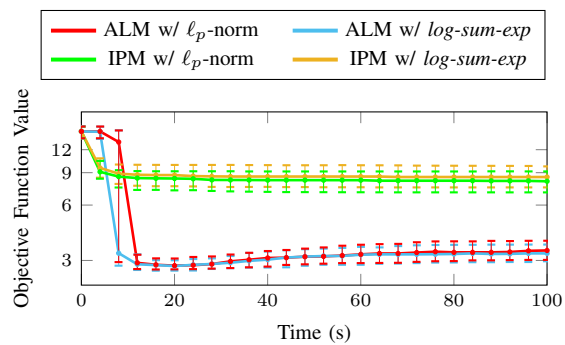
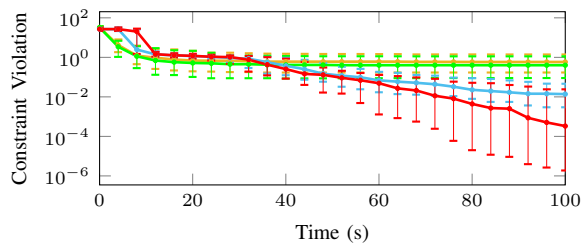 

\begin{figure}[!t] 
    \centering
    \begin{subfigure}{\columnwidth}
    \pgfplotsset{every tick label/.append style={font=\scriptsize}}

\begin{tikzpicture}
\begin{axis}[ymin=0, ymax=0.5,
     xmin=0, xmax=3.5,  
     width=0.9\linewidth,
     height=0.45\linewidth,
        xtick = {0, 1, 2, 3},
        ytick = {0, 0.2, 0.4},
        yticklabels={0,0.2,0.4},
        xlabel={\footnotesize Time (h)},
        ylabel={\footnotesize UAV Remaining ToF (h)}
        ]

\addplot[
draw=none,  fill=gray, fill opacity=0.25] coordinates {
(0.2339,\pgfkeysvalueof{/pgfplots/ymin})
(0.3194,\pgfkeysvalueof{/pgfplots/ymin})
(0.3194,\pgfkeysvalueof{/pgfplots/ymax})
(0.2339,\pgfkeysvalueof{/pgfplots/ymax})
};

\addplot[draw=none,  fill=gray, fill opacity=0.25] coordinates {
(0.591250105088156,\pgfkeysvalueof{/pgfplots/ymin})
(0.694916177413065,\pgfkeysvalueof{/pgfplots/ymin})
(0.694916177413065,\pgfkeysvalueof{/pgfplots/ymax})
(0.591250105088156,\pgfkeysvalueof{/pgfplots/ymax})
};

\addplot[draw=none,  fill=gray, fill opacity=0.25] coordinates {
(0.794841028885891,\pgfkeysvalueof{/pgfplots/ymin})
(0.855139960683371,\pgfkeysvalueof{/pgfplots/ymin})
(0.855139960683371,\pgfkeysvalueof{/pgfplots/ymax})
(0.794841028885891,\pgfkeysvalueof{/pgfplots/ymax})
};

\addplot[draw=none,  fill=gray, fill opacity=0.25] coordinates {
(1.01116401009824,\pgfkeysvalueof{/pgfplots/ymin})
(1.08302742792247,\pgfkeysvalueof{/pgfplots/ymin})
(1.08302742792247,\pgfkeysvalueof{/pgfplots/ymax})
(1.01116401009824,\pgfkeysvalueof{/pgfplots/ymax})
};

\addplot[draw=none,  fill=gray, fill opacity=0.25] coordinates {
(1.19462364830708,\pgfkeysvalueof{/pgfplots/ymin})
(1.46216356425611,\pgfkeysvalueof{/pgfplots/ymin})
(1.46216356425611,\pgfkeysvalueof{/pgfplots/ymax})
(1.19462364830708,\pgfkeysvalueof{/pgfplots/ymax})
};

\addplot[draw=none,  fill=gray, fill opacity=0.25] coordinates {
(1.83511146006658,\pgfkeysvalueof{/pgfplots/ymin})
(1.98667815144218,\pgfkeysvalueof{/pgfplots/ymin})
(1.98667815144218,\pgfkeysvalueof{/pgfplots/ymax})
(1.83511146006658,\pgfkeysvalueof{/pgfplots/ymax})
};

\addplot[draw=none,  fill=gray, fill opacity=0.25] coordinates {
(2.23904367015298,\pgfkeysvalueof{/pgfplots/ymin})
(3.18442101025490,\pgfkeysvalueof{/pgfplots/ymin})
(3.18442101025490,\pgfkeysvalueof{/pgfplots/ymax})
(2.23904367015298,\pgfkeysvalueof{/pgfplots/ymax})
};

\addplot[draw=none,  fill=gray, fill opacity=0.25] coordinates {
(3.35008194849455,\pgfkeysvalueof{/pgfplots/ymin})
(3.37996329245054,\pgfkeysvalueof{/pgfplots/ymin})
(3.37996329245054,\pgfkeysvalueof{/pgfplots/ymax})
(3.35008194849455,\pgfkeysvalueof{/pgfplots/ymax})
};

    \addplot[solid, very thick, color = blue, line join=round]
      table[x=time, y=t_A] {battery_refined.dat};
        
    \addplot[
      only marks,
      mark=asterisk,
      mark options={draw=red, line width=1pt, scale=1.5}
    ]
    coordinates {
    (0.120421207048964,0.279578598685544)
    (0.455554720598500,0.158142803053360)
    (0.744813747031028,0.120690805095128)
    (0.933001818734942,0.0831732936330849)
    (1.13106481631317,0.0647059652814087)
    (1.64484823822156,0.216882157673670)
    (1.69950989421480,0.162220287977078)
    (1.81696883720691,0.0447539794814092)
    (2.11186877240417,0.128199777268160)
    (3.25888100376545,0.325539967758289)
    };

\end{axis}
\end{tikzpicture}
    \caption{History of UAV's remaining time-of-flight (ToF).}
    \label{subfigure: history of UAV battery level}
    \end{subfigure}

    \bigskip

    \begin{subfigure}{\columnwidth}
    \pgfplotsset{every tick label/.append style={font=\scriptsize}}

\begin{tikzpicture}
\begin{axis}[ymin=-0.5, ymax=6.5,
     xmin=0, xmax=3.5,  
     width=0.9\linewidth,
     height=0.45\linewidth,
        xtick = {0, 1, 2, 3},
        ytick = {0, 3.0, 6.0},
        yticklabels={0, 3.0, 6.0},
        xlabel={\footnotesize Time (h)},
        ylabel={\footnotesize UAV-UGV Distance (km)}
        ]

\addplot[
draw=none,  fill=gray, fill opacity=0.25] coordinates {
(0.2339,\pgfkeysvalueof{/pgfplots/ymin})
(0.3194,\pgfkeysvalueof{/pgfplots/ymin})
(0.3194,\pgfkeysvalueof{/pgfplots/ymax})
(0.2339,\pgfkeysvalueof{/pgfplots/ymax})
};

\addplot[draw=none,  fill=gray, fill opacity=0.25] coordinates {
(0.591250105088156,\pgfkeysvalueof{/pgfplots/ymin})
(0.694916177413065,\pgfkeysvalueof{/pgfplots/ymin})
(0.694916177413065,\pgfkeysvalueof{/pgfplots/ymax})
(0.591250105088156,\pgfkeysvalueof{/pgfplots/ymax})
};

\addplot[draw=none,  fill=gray, fill opacity=0.25] coordinates {
(0.794841028885891,\pgfkeysvalueof{/pgfplots/ymin})
(0.855139960683371,\pgfkeysvalueof{/pgfplots/ymin})
(0.855139960683371,\pgfkeysvalueof{/pgfplots/ymax})
(0.794841028885891,\pgfkeysvalueof{/pgfplots/ymax})
};

\addplot[draw=none,  fill=gray, fill opacity=0.25] coordinates {
(1.01116401009824,\pgfkeysvalueof{/pgfplots/ymin})
(1.08302742792247,\pgfkeysvalueof{/pgfplots/ymin})
(1.08302742792247,\pgfkeysvalueof{/pgfplots/ymax})
(1.01116401009824,\pgfkeysvalueof{/pgfplots/ymax})
};

\addplot[draw=none,  fill=gray, fill opacity=0.25] coordinates {
(1.19462364830708,\pgfkeysvalueof{/pgfplots/ymin})
(1.46216356425611,\pgfkeysvalueof{/pgfplots/ymin})
(1.46216356425611,\pgfkeysvalueof{/pgfplots/ymax})
(1.19462364830708,\pgfkeysvalueof{/pgfplots/ymax})
};

\addplot[draw=none,  fill=gray, fill opacity=0.25] coordinates {
(1.83511146006658,\pgfkeysvalueof{/pgfplots/ymin})
(1.98667815144218,\pgfkeysvalueof{/pgfplots/ymin})
(1.98667815144218,\pgfkeysvalueof{/pgfplots/ymax})
(1.83511146006658,\pgfkeysvalueof{/pgfplots/ymax})
};

\addplot[draw=none,  fill=gray, fill opacity=0.25] coordinates {
(2.23904367015298,\pgfkeysvalueof{/pgfplots/ymin})
(3.18442101025490,\pgfkeysvalueof{/pgfplots/ymin})
(3.18442101025490,\pgfkeysvalueof{/pgfplots/ymax})
(2.23904367015298,\pgfkeysvalueof{/pgfplots/ymax})
};

\addplot[draw=none,  fill=gray, fill opacity=0.25] coordinates {
(3.35008194849455,\pgfkeysvalueof{/pgfplots/ymin})
(3.37996329245054,\pgfkeysvalueof{/pgfplots/ymin})
(3.37996329245054,\pgfkeysvalueof{/pgfplots/ymax})
(3.35008194849455,\pgfkeysvalueof{/pgfplots/ymax})
};

    \addplot[solid, very thick, color = blue, line join=round]
      table[x=time, y=d] {distance_history.dat};
        
    \addplot[
      only marks,
      mark=asterisk,
      mark options={draw=red, line width=1pt, scale=1.5}
    ]
    coordinates {
    (0.120421207048964,4.10144691080170)
    (0.455554720598500,4.88546112095232)
    (0.744813747031028,1.79467992036438)
    (0.933001818734942,2.80785659353481)
    (1.13106481631317,1.80313712569945)
    (1.64484823822156,6.02938880885817)
    (1.69950989421480,3.64047364873568)
    (1.81696883720691,0.625609729150954)
    (2.11186877240417,3.86816314910544)
    (3.25888100376545,1.54804918275115)
    };

\end{axis}
\end{tikzpicture}
    \caption{History of UAV-UGV distance (km).}
    \end{subfigure}
\caption{UAV-UGV trajectory computed via NLP (ALM with $\ell_p$-norm) for the problem shown in Fig.~\ref{fig: an_example}. Each red asterisk marks the point along the trajectory where the UAV reaches a task location.}
\label{fig: uav_battery}
\end{figure}
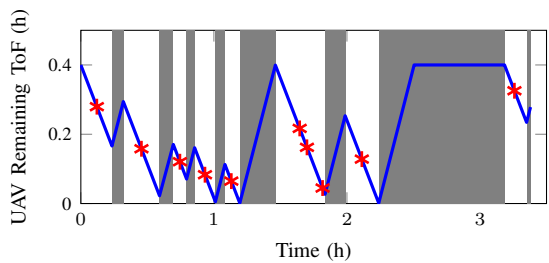
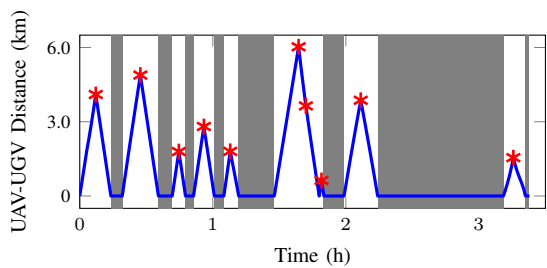

Fig.~\ref{fig: uav_battery} shows the history of the UAV's remaining time of flight as well as the distance between the UAV and the UGV along an optimal trajectory (computed using the ALM with the $\ell_p$-norm approximation) for the problem instance illustrated in Fig.~\ref{fig: an_example}.
The shaded regions indicate periods during which the UAV is charging on the UGV. 
Fig.~\ref{subfigure: history of UAV battery level} shows that the proposed model in \eqref{opt: NLP} allows both partial recharging of the UAV and covering multiple task locations during one discharging cycle of the UAV.

\subsection{Comparison against MINLP}
We compare the performance of the proposed NLP model, solved using the ALM with the $\ell_p$-norm approximation, against the MINLP model. For the ALM, we use the same implementation as the one in Section~\ref{subsec: NLP algorithm}. For MINLP algorithms, we use the open-source solver \texttt{Juniper}~\cite{Juniper} in the \texttt{JuMP}~\cite{JuMP} framework. This solver implements the Branch-and-Bound algorithm with \texttt{IPOPT}~\cite{IPOPT} as the inner NLP solver (with limited-memory Hessian approximation).

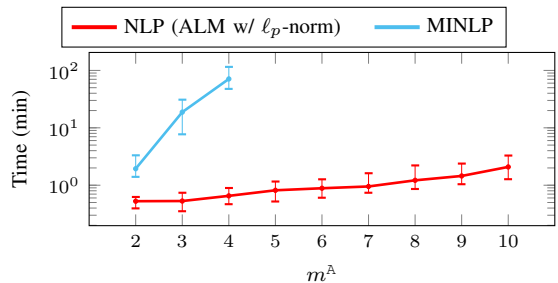
\begin{figure}[!t]
\centering

\begin{tikzpicture}[baseline]
\pgfplotslegendfromname{sharedlegend_1}
\end{tikzpicture}

\pgfplotsset{every tick label/.append style={font=\scriptsize}}
\begin{tikzpicture}

\begin{semilogyaxis}[
    xlabel = \footnotesize \(m^\texttt{A}\), ylabel=\footnotesize Time (min),
    width= 0.9\linewidth,
    height = 0.45\linewidth,
    xmin = 1, xmax = 11,
    xtick = {2,3,4,5,6,7,8,9,10},
    legend to name = sharedlegend_1, 
    legend style={
        legend columns= 2,
        at={(0.5, 0.6)},
        anchor=north west,
        /tikz/column 2/.style={
            column sep=10pt},
        font=\footnotesize
    },
    legend image post style={line width=1.5pt} 
]

\definecolor{NLP_color}{rgb}{1,0,0}
\addplot [mark=none, color=NLP_color,line width=1pt] table[x=aois,y=objMedianVals] {NLP.dat};
\addlegendentry{NLP (ALM w/ \(\ell_p\)-norm)}

\definecolor{MINLP_color}{rgb}{0.301, 0.745, 0.933}
\addplot [mark=none, color=MINLP_color,line width=1pt] table[x=aois,y=objMedianVals] {MINLP.dat};
\addlegendentry{MINLP}

\addplot[
    only marks,
    mark=*,
    mark size=0.8pt,
    color=NLP_color,
    error bars/.cd,
        y dir=both,
        y explicit,
        error bar style={line width=0.1pt},
        error mark=-,  
        error mark options={
            line width=0.8pt,
            mark size=1.5pt,
            rotate=90
        },
]
table[
    x=aois,
    y=objMedianVals,
    y error plus expr=\thisrow{objUpperMedianVals}-\thisrow{objMedianVals},
    y error minus expr=\thisrow{objMedianVals}-\thisrow{objLowerMedianVals},
]{NLP.dat};

\addplot[
    only marks,
    mark=*,
    mark size=0.8pt,
    color=MINLP_color,
    error bars/.cd,
        y dir=both,
        y explicit,
        error bar style={line width=0.1pt},
        error mark=-,  
        error mark options={
            line width=0.8pt,
            mark size=1.5pt,
            rotate=90
        },
]
table[
    x=aois,
    y=objMedianVals,
    y error plus expr=\thisrow{objUpperMedianVals}-\thisrow{objMedianVals},
    y error minus expr=\thisrow{objMedianVals}-\thisrow{objLowerMedianVals},
]{MINLP.dat};


\end{semilogyaxis}
\end{tikzpicture}
\caption{Median computation time versus \(m^{\texttt{A}}\) over 100 Monte Carlo runs with randomly generated aerial task positions for each value of \(m^{\texttt{A}}\), with a maximum computation time of 6 hours. Error bars indicate the 0.25--0.75 quantiles, and the median constraint violation is on the order of \(10^{-6}\).}
\label{fig: NLP Scaling}
\end{figure}

Table~\ref{tab:aoi_compare} and Fig.~\ref{fig: NLP Scaling} compare the performance of NLP and MINLP over a varying number of UAV task locations $m^\texttt{A}$, with a maximum computation time of 6 hours. For each value of $m^\texttt{A}$, we compare NLP and MINLP using 100 problem instances with randomly generated UAV task locations. 
In Table~\ref{tab:aoi_compare}, we report the success rate and the median objective value, median constraint violation, and median computational time. We define the success rate based on whether the solver returns a valid solution within 6 hours of computational time. In addition, we define the constraint violation for both the NLP and MINLP models as the sum of violations of all constraints in~\eqref{eqn: smooth constr} and~\eqref{eqn: disjunctive constr}. Fig.~\ref{fig: NLP Scaling} further illustrates the scalability of the NLP model compared with the MINLP model.


\begin{table}[!t]
\centering
\caption{Performance comparison between the NLP model~(\ref{opt: NLP}), solved using ALM with the $\ell_p$-norm approximation, and the MINLP model~(\ref{opt: MINLP}), with a maximum computation time of 6 hours.}
\label{tab:aoi_compare}
\begin{tabular}{c | c c c c c}
\toprule
\(m^{\texttt{A}}\) & \textbf{Method} & \textbf{Succ. Rate} & \textbf{Obj.} & \textbf{Con. Vio.} & \textbf{Time (min)} \\
\midrule

\multirow{2}{*}{2}
& NLP  & 100\% & 2.4351 & $5.0\texttt{e}{-7}$ & 0.5 \\
& MINLP & 98\% & 2.4352 & $5.7 \texttt{e}{-7}$ & 1.9 \\
\midrule

\multirow{2}{*}{3}
& NLP   & 100\% & 2.4354 & $2.0 \texttt{e}{-6}$ & 0.5 \\
& MINLP & 99\% & 2.4356 & $2.0 \texttt{e}{-6}$ & 18.8 \\
\midrule

\multirow{2}{*}{4}
& NLP   & 100\% & 2.5188 & $1.9 \texttt{e}{-6}$  & 0.7 \\
& MINLP & 87\% & 2.4352 & $1.9 \texttt{e}{-6}$ & 71.0 \\
\bottomrule
\end{tabular}
\end{table}

Overall, the proposed NLP model outperforms the MINLP model as \(m^\texttt{A}\) increases, especially in computation time. Table~\ref{tab:aoi_compare} shows that under similar level of constraint violation, the NLP model consistently outperforms the MINLP model  in terms of success rate and  computational time, with relatively minor increase in the objective function value when \(m^\texttt{A}=4\). Fig.~\ref{fig: NLP Scaling} further shows that how the computation time scales as \(m^\texttt{A}\) increases (under similar level of constraint violation). The proposed NLP model not only improves the computation time by up to two orders of magnitudes when \(m^\texttt{A}\leq 4\), but also ensure the computation time is within minutes when MINLP cannot produce solutions within hours as \(m^\texttt{A}\) increases. These results showcase the computational benefits of the proposed model against MINLP.  

\section{Conclusion}



We presented a nonlinear trajectory optimization model for energy-sharing UAV-UGV systems with multiple UAV and UGV task locations. We based this model on a smoothing approximation of disjunctive constraints, which eliminates the need for integer programming. We demonstrated the proposed model on a one-UAV-one-UGV system. Compared with mixed-integer nonlinear programming, this model reduces the computation time from hours to minutes in numerical simulation. In future work, we plan to extend the current model to multi-UAV-multi-UGV systems and persistent monitoring applications.







\bibliographystyle{IEEEtran}
\bibliography{IEEEabrv,reference}

\end{document}